*Igor V. Bayak*

# A lemma on the minimal surfaces

Stated lemma contains the assertions about isomorphism of exact m-forms and exterior differentials of regular m-maps, of linearly harmonic m-forms and exterior differentials of regular harmonic m-maps, of global minimal $(n-m)$-surfaces and level $(n-m)$-surfaces of regular minimal m-maps.

Let $I = \{1,..,n\}$; $J = \{1,..,m\}$, $1 \leq m < n$; $(J) \equiv (i_1,...,i_j,...,i_m): J \to I$ - is an arbitrary ranked in ascending order arrangement, i.e. such injection that $i_{j-1} < i_j$; $\{(J)\}$ - is set of all ordered injections; E - is n-dimensional Euclidean space; $(\bar{e}^j)_I$ - is basis of E; $(\hat{e}^{(J)})_{\{(J)\}} \equiv (dx_{i_1} \wedge .. dx_{i_j} .. \wedge dx_{i_m})_{\{(J)\}}$ - is basis of $\Lambda^m E^*$;

$$\varphi^m : E \to R^m \Big|_{\sup\{|\Delta x_j|\}_I \to 0} \lim |\varphi_i^m(\bar{x}+\Delta\bar{x}) - \varphi_i^m(\bar{x})| \to 0$$

$$\lim_{\sup\{|\Delta x_j|\}_I \to 0} \frac{\varphi_i^m(\bar{x}) - \varphi_i^m(\bar{x}+\Delta\bar{x})}{\varphi_i^m(\bar{x}-\Delta\bar{x}) - \varphi_i^m(\bar{x})} \to 1 \; (\forall \bar{x} \in E, i \in J) \text{ -}$$

is an arbitrary m-map, where simultaneously are satisfied the conditions of continuity and existence of differential. If $J(\bar{x}) \equiv \left(\frac{\partial \varphi_i^m(\bar{x})}{\partial x_j}\right)_{i \in J, j \in I}$ - is Jacobi matrix of the map $\varphi^m$ in the point $\bar{x}$ and rank $J(\bar{x}) = m \; \forall \bar{x} \in E$ (that is equivalent to the existence of m linear independent gradients in all Euclidean space), then m-map is called by regular in E. Let $\{\varphi^m\}$ - is space of all regular m-maps, then $\varphi_{(n-m)}(\bar{x}_0) \equiv \{\bar{x} \mid \varphi^m(\bar{x}) = \varphi^m(\bar{x}_0)\}$ - is level $(n-m)$-surface corresponding to the regular map $\varphi^m$ and to the point of belonging $\bar{x}_0$; $\varphi_{(n-m)} \equiv \{\varphi_{(n-m)}(\bar{x}_0)\}_{\bar{x}_0 \in E}$ - is $(n-m)$-foliation of E and $\{\varphi_{(n-m)}\}$ - is space of $(n-m)$-foliations of E. In addition we define that $d\varphi^m(\bar{x}) \equiv \sum_{\{(J)\}} \det\left(\frac{\partial \varphi_i^m(\bar{x})}{\partial x_j}\right)_{j \in (J)}^{i \in J} \hat{e}^{(J)}$ - is exterior differential of $\varphi^m$ in the point $\bar{x}$; $d\varphi^m \equiv \{d\varphi^m(\bar{x})\}_{\bar{x} \in E}$ - is exterior differential of $\varphi^m$ in E; $\{d\varphi^m\}$ - is space of exterior differentials in E. Let $H^i(\bar{x}) \equiv \left(\frac{\partial^2 \varphi_i^m(\bar{x})}{\partial x_j \partial x_k}\right)_{j,k \in I}$, where $i \in J$, - are Hesse matrices of $\varphi^m$ in the point $\bar{x}$, then $\{\eta^m\} \equiv \{\varphi^m \mid \text{Sp}H^i = 0 \; (\forall i \in J)\}$ - is space of all regular harmonic m-maps. Let also $H^i_{(J)}(\bar{x}) \equiv \left(\frac{\partial^2 \varphi_i^m(\bar{x})}{\partial x_j \partial x_k}\right)_{j,k \in I \setminus (J)}$, where $i \in J, (J) \in \{(J)\}$, - are sectional Hesse matrices of $\varphi^m$ in the point $\bar{x}$, then $\{\mu^m\} \equiv \{\varphi^m \mid \text{Sp}H^i_{(J)} = 0 \; (\forall i \in J, (J) \in \{(J)\})\}$ - is space of all regular minimal m-maps, in particular, $\{\mu^1\} \equiv \{\varphi^1 \mid \mathrm{SpH}_i^1 = 0 \, (\forall \, i \in I)\}$ - is space of all regular minimal 1-maps. Further let $\{\hat{\varphi}^m\}$ - is space of all exact m-forms and $\{\hat{\eta}^m\}$ - is space of all linearly harmonic m-forms, i.e. such m-forms that vanish by linear operator $(d + \delta)$, and then hold the affirmation.

**Lemma:**

1. $\{d\varphi^m\} \approx \{\hat{\varphi}^m\}$
2. $\{d\eta^m\} \approx \{\hat{\eta}^m\}$
3. Space of foliations $\{\mu_{(n-m)}\}$ is isomorphic to space of foliations of E by global minimal $(n-m)$-surfaces, i.e. the minimality condition of regular m-map serves as the differential criterion of extremality (minimality) of integral functional of volume of surfaces under local variations.